\documentclass[12pt]{article}

\usepackage[english]{babel}

\usepackage{amsfonts}
\usepackage{amssymb}
\usepackage{amsmath}
\usepackage{amsthm}
\usepackage{epsf}
\usepackage{float}

\usepackage{cite}

\usepackage{fullpage}

\newtheorem{theorem}{Theorem}
\newtheorem{prp}[theorem]{Proposition}

\theoremstyle{definition}

\newcommand{\bbR}{{\mathbb R}}
\newcommand{\bbN}{{\mathbb N}}

\newcommand{\lr}[1]{\!\left(\!#1\!\right)}

\newcommand{\cmin}{c^{\rm min}}
\newcommand{\cmax}{c^{\rm max}}

\newcommand{\dx}{dx}

\title{Error estimates for Riemann sums of some singular functions}
\author{Pavel  Gurevich\footnote{Free University of Berlin, RUDN University, email: gurevich@math.fu-berlin.de}, \  
        Sergey  Tikhomirov\footnote{Saint-Petersburg State Univeristy; email: s.tikhomirov@spbu.ru}}

\begin{document}

\maketitle

\begin{abstract}
In this short note, we obtain error estimates for Riemann sums of some singular functions.
\end{abstract}

Let $N \in \bbN$, and let $z_0 =-1$ or $z_0 = 0$.
For a function $F_1 :[z_0, 1] \to \bbR$ denote the error estimate of the Riemann sum of the integral $\int_{z_0}^1 F_1(x) dx$ by
$$
R_n := \int_{z_0}^1 F_1(x) dx - \left( \frac{1}{2n}F_1(z_0) + \frac{1}{2n}F_1(1) + \sum_{ k = z_0 n + 1}^{n-1} \frac{1}{n}F_1\lr{\dfrac{k}{n}}\right).
$$

We prove the following 3 propositions, which are used, e.g., in~\cite{GT}.

\begin{prp}\label{lemIntSqrt1}
Assume that a function $F_1(x)$ can be represented as
$$
F_1(x) = c_1 (1-x)^{1/2}+ c_2 (1-x)^{3/2} + \tilde{F}_1(x),
$$
where $c_1, c_2 \in \bbR$ and $\tilde{F}_1 \in C^2[z_0, 1]$. Then

\begin{enumerate}
\item
There exists $L_1 = L_1(F_1, N) > 0$ such that
$$
\left|R_n\right| \leq L_1\frac{1}{n^{3/2}}, \quad n \geq N.
$$
\item If, additionally, $\tilde{F}_1(x) = c_3 (1-x)^{5/2} + c_4 (1-x)^{7/2} + \bar{F}_1(x)$, where $c_3, c_4 \in \bbR$ and $\bar{F} \in C^4[z_0, 1]$, then there exists $\bar{L}_1 = \bar{L}_1(F_1, N) > 0$ such that
$$
\left|(n+1)^2 R_{n+1} - n^2 R_n\right| \leq \bar{L}_1\frac{1}{n^{1/2}} \quad n \geq N.
$$
\end{enumerate}
\end{prp}

\begin{prp}\label{lemIntSqrt2}
Assume that a function $F_2(x)$ can be represented as
\begin{equation}\label{P7s}
F_2(x) = c_1 (1-x)^{-1/2}+ c_2 (1-x)^{1/2} + \tilde{F}_2(x),
\end{equation}
where $c_1 > 0$, $c_2 \in \bbR$, and $\tilde{F}_2 \in C^1[z_0, 1]$.
Then there exists $L_2 = L_2(F_2, N) > 0$, $L_2^* = L_2^*(F_2, N) > c_1$ and $l_2 = l_2(F_2, N)$ such that
\begin{equation}\label{eqL2}
L_2^*\frac{1}{n^{1/2}} - l_2 \frac{1}{n} \leq \int_{z_0}^1 F_2(x) dx - \sum_{k = z_0n}^{n-1} \frac{1}{n}F_2\lr{\dfrac{k}{n}} \leq
L_2\frac{1}{n^{1/2}}, \quad n \geq N.
\end{equation}
In particular, $\left|\sum \limits_{k = z_0 n}^{n-1} \frac{1}{n}F_2\lr{\frac{k}{n}}\right|$ is bounded.
\end{prp}

\begin{prp}\label{propL3}
Assume that $F_3 \in C^0[-1, 1)$ and $|(1-x)^{3/2}F_3(x)|$ is bounded on $[-1, 1)$. Then there exists $L_3 = L_3(F_3, N) > 0$ such that
$$
\left|\sum_{|k|\leq n-1} \frac{1}{n}F_3\lr{\frac{k}{n}}\right| \leq L_3{n^{1/2}}, \quad n \geq N.
$$
\end{prp}


Here we provide the proofs of item 2 of Proposition \ref{lemIntSqrt1} and Proposition \ref{lemIntSqrt2} for the case $z_0 = -1$. Item 1 of Proposition \ref{lemIntSqrt1}  and the case $z_0 = 0$ can be proved similarly. Proposition~\ref{propL3} is straightforward.

\subsection*{Proof of item 2 of Proposition \ref{lemIntSqrt1} }

Consider the function $\Phi(x): = F_1(x) + F_1(-x)$.
Then
$$
n^2 R_n = n^2 I_{F1}  - \left( \frac{n}{2}\Phi(0) + \frac{n}{2}\Phi(1) + \sum_{k=1}^{n-1} n \Phi\lr{\frac{k}{n}} \right),
$$
where $I_{F1} := \int_{-1}^1 F_1(x)dx$. Therefore,
\begin{align}
(n+1)^2R_{n+1} - n^2 R_n
& = (2n+1)I_{F1} -
\left(
\frac{\Phi(0) + \Phi(1)}{2} +
\sum_{k = 1}^{n} (n+1) \Phi\lr{\frac{k}{n+1}} -
\sum_{k = 1}^{n-1} n \Phi\lr{\frac{k}{n}}
\right) \notag \\
& = \Sigma_{1, n} + \Sigma_{2, n}, \label{P13s}
\end{align}
where
\begin{align*}
\Sigma_{1, n} & := (n+1)I_{F1} - \left( \frac{\Phi(0) + \Phi(1)}{2} + \sum_{k = 1}^n \Phi\lr{\frac{k}{n+1}} \right),\\
\Sigma_{2, n} & := n\lr{I_{F1} - \Phi(0) - \Sigma_{\Phi, n} }, \quad
\Sigma_{\Phi, n}  := \sum_{k = 1}^{n} \left( \Phi\lr{\frac{k}{n+1}} - \Phi\lr{\frac{k-1}{n}} \right).
\end{align*}
Note that $\Sigma_{1, n} = (n+1)R_{n+1}$. Hence, Proposition \ref{lemIntSqrt1} (item 1) implies that, for some $L_1 > 0$,
\begin{equation}\label{P14ss}
|\Sigma_{1, n}| \leq L_1 \frac{1}{\sqrt{n+1}}.
\end{equation}

Estimates of $\Sigma_{2, n}$ are more subtle. First, note that
$\frac{k}{n+1} - \frac{k-1}{n} = \frac{n+1-k}{n(n+1)}$
and, hence,
\begin{align*}
\Phi\lr{\frac{k}{n+1}} - \Phi\lr{\frac{k-1}{n}} & =  \frac{n+1-k}{n(n+1)}\Phi'\lr{\frac{k}{n+1}} \\ & \;\;\;\; - \frac{1}{2}\left(\frac{n+1-k}{n(n+1)}\right)^2 \Phi''\lr{\frac{k}{n+1}} +
\frac{1}{6}\left(\frac{n+1-k}{n(n+1)}\right)^3 \Phi'''(\xi_k),
\end{align*}
where $\xi_k \in \left[\frac{k-1}{n}, \frac{k}{n+1}\right]$. Therefore,
\begin{align}
\Sigma_{\Phi, n} & =  \frac{n+1}{n}\sum_{k = 1}^{n}\frac{1}{n+1}\frac{n+1-k}{n+1}\Phi'\lr{\frac{k}{n+1}}   -\frac{1}{2}\frac{n+1}{n^2}\sum_{k = 1}^{n}\frac{1}{n+1}\lr{\frac{n+1-k}{n+1}}^2\Phi''\lr{\frac{k}{n+1}} \notag \\
& \;\;\;\;  + \frac{1}{6}\frac{n+1}{n^3}\sum_{k = 1}^{n}\frac{1}{n+1}\lr{\frac{n+1-k}{n+1}}^3\Phi'''\lr{\xi_k}. \label{P14s}
\end{align}
Set
\begingroup
\allowdisplaybreaks
$$
\Phi_1(x) := (1-x)\Phi'(x), \quad \Phi_2(x) := (1-x)^2\Phi''(x), \quad \Phi_3(x) := (1-x)^3\Phi'''(x),
$$
\begin{equation}\label{P'16s}
I_{\Phi 1}:= \int_{0}^1 \Phi_1(x) dx, \quad I_{\Phi 2}:= \int_{0}^1 \Phi_2(x) dx,
\end{equation}
\begin{align}
\Sigma_{\Phi 1, n} & := I_{\Phi 1} - \sum_{k = 1}^{n}\frac{1}{n+1}\frac{n+1-k}{n+1}\Phi'\lr{\frac{k}{n+1}} - \frac{1}{2(n+1)}\Phi'(0), \label{eqf1}\\
\Sigma_{\Phi 2, n} & := I_{\Phi 2} - \sum_{k = 1}^{n}\frac{1}{n+1}\lr{\frac{n+1-k}{n+1}}^2\Phi''\lr{\frac{k}{n+1}} - \frac{1}{2(n+1)}\Phi''(0), \label{eqf2}\\
\Sigma_{\Phi 3, n} & := \sum_{k = 1}^{n}\frac{1}{n+1}\lr{\frac{n+1-k}{n+1}}^3\Phi'''\lr{\xi_k} \label{eqf3}.
\end{align}
\endgroup
Since $\xi_k \leq \dfrac{k}{n+1}$, we have
\begin{equation}\label{P15s}
\left|\lr{\frac{n+1-k}{n+1}}^3\Phi'''\lr{\xi_k}\right| \leq |\Phi_3(\xi_k)|.
\end{equation}


Proposition \ref{lemIntSqrt1} (item 1) and \eqref{P15s} imply that
\begin{equation}\label{eqKf}
|\Sigma_{\Phi 1, n}| \leq K_{\Phi 1}\frac{1}{(n+1)^{3/2}}, \quad  |\Sigma_{\Phi 2, n}| \leq K_{\Phi 2}\frac{1}{(n+1)^{3/2}}, \quad |\Sigma_{\Phi 3, n}| \leq K_{\Phi 3},
\end{equation}
$$
K_{\Phi 1} := L_1(\Phi_1, N), \quad K_{\Phi 2} := L_1(\Phi_2, N), \quad K_{\Phi 3} := \sup_{x \in (0, 1)} |\Phi_3(x)|,
$$

In this notation, \eqref{P14s} takes the form
\begin{equation*}\notag
\Sigma_{\Phi, n}  =  \frac{n+1}{n}\lr{I_{\Phi 1} - \frac{\Phi'(0)}{2(n+1)} - \Sigma_{\Phi 1, n}}
  -\frac{n+1}{2n^2}\lr{I_{\Phi 2}- \frac{\Phi''(0)}{2(n+1)} - \Sigma_{\Phi 2, n}}  + \frac{n+1}{6n^3}\Sigma_{\Phi 3, n}.
\end{equation*}
Finally, taking into account \eqref{P14ss},
we conclude that
\begin{align*}
\Sigma_{2, n} & = nI_{F1} - (n+1)I_{\Phi 1}+ \frac{1}{2}\frac{n+1}{n}I_{\Phi 2}  - \lr{n\Phi(0) - \frac{1}{2}\Phi'(0) + \frac{1}{4n} \Phi''(0)} \\
& \;\;\;\; + \lr{(n+1)\Sigma_{\Phi 1, n} - \frac{1}{2}\frac{n+1}{n}\Sigma_{\Phi 2, n} - \frac{1}{6}\frac{n+1}{n^2}\Sigma_{\Phi 3, n}}.
\end{align*}
Using the relations $I_{F1} = I_{\Phi 1} + \Phi(0)$ and $I_{\Phi 2} = 2I_{\Phi 1} - \Phi'(0)$, we obtain
\begin{align*}
\Sigma_{2, n} & = nI_{\Phi 1} + n\Phi(0) -(n+1)I_{\Phi 1}+ \frac{n+1}{n}I_{\Phi 1} - \frac{1}{2}\frac{n+1}{n}\Phi'(0) \\
 & \;\;\;\;  -\lr{n\Phi(0) - \frac{1}{2}\Phi'(0) + \frac{1}{4n}\Phi''(0)}  + \lr{(n+1)\Sigma_{\Phi 1, n} - \frac{1}{2}\frac{n+1}{n}\Sigma_{\Phi 2, n} - \frac{1}{6}\frac{n+1}{n^2}\Sigma_{\Phi 3, n}} \\
 & = \frac{1}{n}I_{\Phi 1} - \frac{1}{4n}\Phi''(0) - \frac{1}{2n}\Phi'(0) + \lr{(n+1)\Sigma_{\Phi 1, n} - \frac{1}{2}\frac{n+1}{n}\Sigma_{\Phi 2, n} - \frac{1}{6}\frac{n+1}{n^2}\Sigma_{\Phi 3, n}}.
\end{align*}
The latter equality and estimates \eqref{eqKf} yield
\begin{align*}
|\Sigma_{2, n}| &  \leq \frac{1}{n} \lr{|I_{\Phi 1}| + \frac{|\Phi''(0)|}{4} + \frac{|\Phi'(0)|}{2} + \frac{K_{\Phi 3}}{3} } + K_{\Phi 1}\frac{1}{n^{1/2}} + \frac{K_{\Phi 2}}{2}\frac{1}{n^{3/2}} \\
& \leq \frac{1}{n^{1/2}}\lr{K_{\Phi 1} + \lr{|I_{\Phi 1}| + \frac{|\Phi''(0)|}{4} + \frac{|\Phi'(0)|}{2} + \frac{K_{\Phi 3}}{3}}\frac{1}{N^{1/2}} + \frac{K_{\Phi 2}}{2}\frac{1}{N}}.
\end{align*}
Combining this with \eqref{P13s} and \eqref{P14ss}, we have
\begin{equation}\notag
|(n+1)^2 R_{n+1} - n^2 R_n| \leq \bar{L}_1 \frac{1}{\sqrt{n}},\end{equation}
where

\begin{equation}\notag
\bar{L}_1 := L_1 + K_{\Phi 1} + \lr{|I_{\Phi 1}| + \frac{|\Phi''(0)|}{4} + \frac{|\Phi'(0)|}{2} + \frac{K_{\Phi 3}}{3}}\frac{1}{N^{1/2}} + K_{\Phi 2}\frac{1}{N}.
\end{equation}

\subsection*{Proof of Proposition \ref{lemIntSqrt2}}

Fix $\cmax_1 > c_1$ and $\cmin_1 \in (\frac{2}{3\sqrt{2}-2}c_1, c_1)$.
Let us choose $\delta \in (0, 1)$  such that, for $x \in [1-\delta, 1)$,
\begin{equation}\label{eqIntSqrt01}
\frac{\cmin_1}{(1-x)^{1/2}}  \leq F_2(x)  \leq \frac{\cmax_1}{(1-x)^{1/2}}, \quad \quad
\frac{1}{2}\frac{\cmin_1}{(1-x)^{3/2}}   \leq F'_2(x) \leq \frac{1}{2}\frac{\cmax_1}{(1-x)^{3/2}}.
\end{equation}

Using the representation \eqref{P7s} and \eqref{eqIntSqrt01}, we choose $C = C(F) > 0$ such that
\begin{equation}\label{eqIntSqrtBV}
\int_{-1}^{\beta} |F'_2(x)| \dx \leq F(\beta) + C, \quad \beta \in [1- \delta, 1).
\end{equation}

Finally, we assume that (increasing $L_2$ and $l_2$ if necessarily)
\begin{equation}\label{P33s}
n \geq \frac{2}{\delta}.
\end{equation}

Integrating by parts, we have
\begingroup
\begin{multline*}
\int_{-1}^1 F_2(x) \dx
 = \sum_{k = -n}^{n-1} \int_{k/n}^{(k+1)/n} F_2(x) \dx \\
 = \sum_{k = -n}^{n-1} \lr{\left.\lr{x-\frac{k+1}{n}}F_2(x) \right|_{k/n}^{(k+1)/n} - \int_{k/n}^{(k+1)/n} \lr{x-\frac{k+1}{n}}F'_2(x) \dx}
 = \sum_{k = -n}^{n-1} \frac{1}{n} F_2\lr{\frac{k}{n}} + \Sigma_{1, n},
\end{multline*}
\endgroup
where
$$
\Sigma_{1, n} := \sum_{k = -n}^{n-1} \int_{k/n}^{(k+1)/n} \lr{\frac{k+1}{n}-x}F'_2(x) \dx.
$$
Inequalities \eqref{eqIntSqrt01}--\eqref{P33s} imply that
\begin{align}
|\Sigma_{1, n}|
& = \left|\int_{1-1/n}^{1} (1-x)F'_2(x) + \sum_{k = -n}^{n-2} \int_{k/n}^{(k+1)/n} \lr{\frac{k+1}{n}-x}F'_2(x) \dx \right| \notag \\
& \leq \int_{1-1/n}^{1} \frac{1}{2}\frac{\cmax_1}{(1-x)^{1/2}} + \sum_{k = -n}^{n-2} \frac{1}{n} \int_{k/n}^{(k+1)/n} |F'_2(x)| \dx
 \leq \frac{\cmax_1}{\sqrt{n}} + \frac{1}{n} \int_{-1}^{1-1/n} |F'_2(x)| \dx. \label{P33ss}
\end{align}
Since $1-1/n > 1 - \delta$ due to \eqref{P33s}, we conclude from \eqref{eqIntSqrt01}, \eqref{eqIntSqrtBV} and \eqref{P33ss} that
\begin{equation*}
|\Sigma_{1, n}|  \leq \frac{\cmax_1}{\sqrt{n}} + \frac{1}{n}\lr{F_2\lr{1-\frac{1}{n}} + C} \leq 2\cmax_1\frac{1}{\sqrt{n}} + C\frac{1}{n}.
\end{equation*}
Taking $L_2 := 2\cmax_1 + \frac{C}{\sqrt{N}}$, we obtain the second inequality in \eqref{eqL2}.

Set $k_0 := [(1-\delta) n]$ ($[\cdot]$ is the integer part of a real number). We represent $\Sigma_{1, n}$ as follows:
\begin{equation}\label{eqSigmas}
\Sigma_{1, n} = \Sigma_{2, n} + \Sigma_{3, n},
\end{equation}
where
\begin{align}
\notag \Sigma_{2, n}
& := \lr{\sum_{k = -n}^{k_0-1}  \int_{k/n}^{(k+1)/n} \lr{\frac{k+1}{n}-x}F'_2(x) \dx + \int_{k_0}^{1-\delta} \lr{\frac{k_0+1}{n}-x}F'_2(x) \dx}, \\
\label{eqIntApproxS3} \Sigma_{3, n}
& := \lr{ \int_{1-\delta}^{(k_0+1)/n} \lr{\frac{k_0+1}{n}-x}F'_2(x) \dx + \sum_{k = k_0+ 1}^{n-1} \int_{k/n}^{(k+1)/n} \lr{\frac{k+1}{n}-x}F'_2(x)}.
\end{align}

Below we estimate $\Sigma_{2, n}$ and $\Sigma_{3, n}$ separately.
Inequality \eqref{eqIntSqrtBV} implies that
\begin{equation}\label{eqSigma2IA}
\Sigma_{2, n} \geq -\frac{1}{n}\lr{\int_{-1}^{1-\delta} |F'_2(x)| \dx} \geq -\frac{1}{n}(F_2(1-\delta) + C).
\end{equation}

Due to \eqref{eqIntSqrt01} all the terms in expression $\Sigma_{3, n}$ are positive. Inequality \eqref{P33s} implies that the last sum in formula \eqref{eqIntApproxS3} contains at least two terms and inequality \eqref{eqIntSqrt01} implies that
\begin{align}
\Sigma_{3, n} & \geq \int_{1-2/n}^{1-1/n}\lr{1-\frac{1}{n}-x}F'_2(x) \dx + \int_{1-1/n}^{1}(1-x)F'_2(x) \dx \notag \\
         & \geq \int_{1-2/n}^{1-1/n}\frac{1}{2}\lr{1-\frac{1}{n}-x}\frac{\cmin_1}{(1-x)^{3/2}} \dx + \int_{1-1/n}^{1}\frac{1}{2}\frac{\cmin_1}{(1-x)^{1/2}} \dx \notag\\
         & = \int_{1-2/n}^{1}\frac{1}{2}\frac{\cmin_1}{(1-x)^{1/2}} \dx - \frac{1}{2n}\int_{1-2/n}^{1-1/n}\frac{\cmin_1}{(1-x)^{3/2}} = \frac{3\sqrt{2}-2}{2}\cmin_1\frac{1}{\sqrt{n}} \label{eqSigma3IA}.
\end{align}

Relations \eqref{eqSigmas}, \eqref{eqSigma2IA}, and \eqref{eqSigma3IA} imply that
\begin{equation*}
\Sigma_{1, n} \geq L_2^*\frac{1}{\sqrt{n}} - l_2\frac{1}{n},
\end{equation*}
where $L_2^* := \frac{3\sqrt{2}-2}{2}\cmin_1 > c_1$ and $l_2 := F_2(1-\delta) + C$.


\begin{thebibliography}{99}
\bibitem{GT} Pavel Gurevich, Sergey Tikhomirov, Spatially discrete reaction-diffusion equations with discontinuous hysteresis, 	arXiv:1504.02385
\end{thebibliography}
\end{document}